\documentclass[english]{ourlematema}
\usepackage{amsmath, amsthm, amssymb, hyperref, color}
\usepackage{MnSymbol} 
\usepackage{graphicx}
\usepackage{caption}
\usepackage[all]{xypic}
\usepackage{verbatim}
\usepackage{chemfig, chemnum}
\usepackage{tikz}
\usepackage[linesnumbered,lined,commentsnumbered]{algorithm2e}
\usepackage{caption}
\usepackage{subcaption}
\usepackage{makecell}
\usepackage{array}

 \numberwithin{theorem}{section}

\newtheorem{question}{Question}
\theoremstyle{definition}

\newcommand{\PP}{\mathbb{P}}
\newcommand{\RR}{\mathbb{R}}
\newcommand{\QQ}{\mathbb{Q}}
\newcommand{\CC}{\mathbb{C}}
\newcommand{\ZZ}{\mathbb{Z}}

\newcommand{\bernd}[1]{\strut{\color{blue} $\clubsuit$}%
  \marginpar{\color{blue}\footnotesize\raggedright  \textbf{Bernd:} #1}}
\newcommand{\emre}[1]{\strut{\color{red} $\clubsuit$}%
  \marginpar{\color{red}\footnotesize\raggedright  \textbf{Emre:} #1}}
  \newcommand{\marta}[1]{\strut{\color{green} $\clubsuit$}%
    \marginpar{\color{green}\footnotesize\raggedright  \textbf{Marta:} #1}}

\renewcommand{\emre}[1]{}
\renewcommand{\bernd}[1]{}
\renewcommand{\marta}[1]{}
 
 \title{Twenty-Seven Questions about the Cubic Surface}
 
 \author{Kristian Ranestad}
 \address{%
Department of Mathematics,
University of Oslo \\
\email{ranestad@math.uio.no}
}
\author{Bernd Sturmfels}
\address{%
MPI for Mathematics in the Sciences, Leipzig \\
\email{bernd@mis.mpg.edu }
}

 \date{2019/10/09}

\begin{document}
\maketitle
\begin{abstract}
\noindent
We present a collection of research questions on cubic surfaces
in $3$-space.  These questions inspired a collection
of papers to be published 
in a special issue of the journal {\em Le Matematiche}.
This article serves as the introduction to that issue.
The number of questions is meant to match the number of lines on a cubic surface.
We end with a list of problems that are open.
\end{abstract}

\section{Introduction}

One of most prominent results in classical algebraic geometry,
derived two centuries ago by Cayley and Salmon, states
that every smooth cubic surface in complex projective $3$-space $\PP^3$
contains $27$ straight lines.
This theorem has inspired generations
of mathematicians, and it is still of importance today.
The advance of tropical geometry and computational methods,
and the strong current interest in applications of algebraic geometry,
led us to revisit the cubic~surface.

Section 3 of this article gives a brief exposition of the
history of that classical subject and how it developed.
We offer a number of references that can serve
as first entry points for students of algebraic geometry who
wish to learn more.

\smallskip

In December 2018, the second author compiled the
list of $27$ question. His original text was slightly edited
and it now constitutes our Section 2 below.
These questions were circulated,
and they were studied by various groups
of young mathematicians, in particular in Leipzig and Berlin.
In May  2019, a one-day seminar on cubic surfaces
was held in Oslo. Different teams worked on the questions and
they made excellent progress. It was then decided to
call for a special issue of {\em Le Matematiche}, with 
submission deadline in September 2019.

This call resulted in the 14 articles listed as \cite{BrGe}--\cite{SS}.
In Section 4 we give an overview of this collection, and we
briefly discuss the contribution of each article.
For each article, we identify which of the $27$ question it refers to.
We also highlight ten key problems that remain open.
These make it clear that there is still
a lot of research to be done, and our readers are invited to join the effort.

\section{Questions}

The text that follows 
in this section was written and circulated by the second author in December 2018.
 The list of $27$ questions was conceived  as a working document that evolves over time.
 It was initially not meant to be published. It simply
  represents what the second author wanted to know, but did not know back then.
The current status of the $27$ questions will be our topic in Section~4.

\smallskip

In spite of two centuries of research on cubic surfaces, it
appears that there are still many unresolved questions, especially when it
comes to computational, tropical and applied aspects.
Please feel free to circulate this text.
It is aimed at stimulating further work
on  cubic surfaces, or equivalently, on symmetric $4 \times 4 \times 4 $ tensors.
Your feedback and comments will be greatly appreciated.

\medskip

A cubic surface in $\PP^3$ is the zero set of a homogeneous polynomial
\begin{equation}
\label{eq:cubicf}
 \begin{matrix} f &=& c_1  x^3 + c_2  y^3 + c_3  z^3 + c_4  w^3 + c_5  x^2 y 
+ c_6  x^2 z + c_7  x^2w + c_8  xy^2     \\ & & + \,c_9  y^2 z  + c_{10}  y^2w 
+ c_{11}  xz^2   + c_{12}  yz^2 + c_{13}  z^2w + c_{14}  x w^2 \\ & & + \,c_{15}  yw^2  
  + c_{16}  zw^2 + c_{17} xyz + c_{18}  xyw + c_{19}  xzw + c_{20}  yzw.
  \end{matrix}
\end{equation}
We work over a field $K$ of characteristic $0$, such as $\QQ$, $\QQ_p$,
$\QQ(t)$, $\RR$, $\CC$, $\overline{\QQ(t)}$, or $\CC \{\!\{ t \} \!\}$.
The  $15$-dimensional group ${\rm PGL}(4)$ acts naturally on
 the projective space $\PP^{19} = \PP( {\rm Sym}_3(K^4))$ whose coordinates
 are $(c_1:c_2:\cdots : c_{20})$. Our first question concerns the orbits of
 that action.
Here the point of departure would be Kazarnovskii's general formula, found in
\cite{DK, Kaz}, for degrees of orbits.

\begin{question}
Given a generic homogeneous cubic $f$ in $x,y,z,w$, what can we say about the orbit closure
$\,\overline{{\rm PGL}(4) \cdot f}$? What is the degree of this variety in $\PP^{19}$?
Can we determine some of its defining polynomial equations?
\end{question}

The geometric invariant theory of cubic surfaces is well understood.
In his 1861 article \cite{Sal},
Salmon found  six fundamental invariants. Their
degrees are $8$, $16$, $24$, $32$, $40$ and $100$. The
square of the last one is a polynomial in the first five.
Over a century later,
Beklemishev \cite{Bek} proved that Salmon's list is complete.

\begin{question}
How to evaluate the six invariants? \ Same question for their tropicalizations.
   (See the textbook \cite{MS} for basics on tropical geometry). 
\end{question}

One obvious invariant is the discriminant of $f$. This is a polynomial of degree
$32$ in the coefficients $c_1,\ldots,c_{20}$.
Edge \cite{edge} corrects a formula written in 
fundamental invariants due to Salmon, apparently also repeated by Clebsch.

The next question assumes familiarity with the combinatorial theory of
discriminants that was developed by
Gel'fand, Kapranov and Zelevinsky  \cite{GKZ}.

\begin{question}
How many monomials does the discriminant have? How many vertices
does its Newton polytope have, i.e.~how many 
D-equivalence classes of regular triangulations are there?
\end{question}

The discriminant of the cubic $f$ equals the resultant of the quadrics
$\frac{\partial f}{\partial x}$, $\frac{\partial f}{\partial y}$,  $\frac{\partial f}{\partial z}$,
$\frac{\partial f}{\partial w}$. We therefore seek formulas for the resultant of four quadrics
in $x,y,z,w$.

\begin{question}
Can we write the resultant of four quadratic surfaces in $\PP^3$
as the determinant of an $8\times 8$-matrix whose entries are
linear forms in the brackets? This resultant 
is the Chow form of the Veronese embedding of $\PP^3$ into $\PP^9$.
Such a formula would be derived from a nice Ulrich sheaf on that Veronese threefold.
\end{question}

Assuming the answer to Question~4 to be affirmative,
we specialize to get an $8 \times 8$-matrix
whose entries are quartics in $c_1,\ldots,c_{20}$ and whose
determinant equals the discriminant of $f$. Note that the discriminant has degree
$32$ in the $c_i$.

\begin{question}
Which varieties in $\PP^{19}$ arise by imposing
rank conditions on the $8 {\times} 8$-matrix in Question 4?
\end{question}

In his 1899 article \cite{Nan}, E.J.~Nansen writes the above resultant as
the determinant of a $20 \times 20$-matrix.
We can again specialize this to a matrix whose
determinant is the discriminant of~$f$.

\begin{question}
Which varieties in $\PP^{19}$ arise by imposing
rank conditions on this $20 {\times} 20$-matrix?
\end{question}

It seems reasonable to surmise that the
loci in Questions 5 and 6 are cubic surfaces with 
prescribed types of singularities. The simplest
scenario is the occurrence of simple nodes as the only singularities. Then these are at most $4$.  

\begin{question}
For $k=2,3,4$, the variety of $k$-nodal cubics is irreducible
of codimension $k$ in $\PP^{19}$. Using his numerical software \cite{BT},
Sascha Timme
computed that the degrees of these varieties are
$280$, $800$ and $305$ respectively.
Later we learned that these degrees, and many more, 
had been found  by Vainsencher~\cite{Vai}.
Can we find explicit low-degree polynomials that vanish on these
varieties?
\end{question}

\begin{question}
Can we find $17$ real points in $\PP^3$ such that
all $280$ of the $2$-nodal cubics through these points are real?
Can we find $16$ real points in $\PP^3$ such that
all $800$ of the $3$-nodal cubics through these points are real?
Also, are there configurations such that no such cubic is real?
\end{question}

The $4$-nodal cubics are Cayley symmetroids.
These arise in convex optimization, as boundaries of feasible regions in
semidefinite programming.

\begin{question}
Can we find  $15$ real points in $\PP^3$ that lie on $305$ 
real Cayley symmetroids?
\end{question}

The following question arose from a conversation with
Hannah Markwig.

\begin{question}
Can the numbers $280$, $800$ and $305$ be derived tropically?
\end{question}

The next question refers to the classical construction of cubic surfaces
by blowing up the projective plane $\PP^2$ in six points.

\begin{question}
How to construct six points  with integer coordinates in $\PP^2$,
and a basis for the space of cubics  through these points,
such that the resulting cubic surface in $\PP^3$ has
a smooth tropical surface for its $2$-adic tropicalization?
Which unimodular triangulations of the Newton polytope $3 \Delta_3$ of a
dense cubic  $f$ arise?
\end{question}

Up to symmetry, there are 14373645 unimodular triangulations of 
the triple tetrahedron $3 \Delta_3$. This number was reported recently
in~\cite[Theorem 3.1]{JPS}.

\begin{question}
Given a cubic surface over a valued field $K$, how to decide whether
its tropicalization is smooth after some linear change of coordinates?
How to search ${\rm PGL}(4,K)$?
\end{question}

Salmon's invariant of degree $100$ vanishes precisely when
the cubic surface has an Eckardt point, that is, a point  common
to three of its $27$ lines. This invariant deserves further study.

\begin{question}
What is the singular locus of the  Eckardt
hypersurface of degree $100$ in $\PP^{19}$?
\end{question}

A homogeneous cubic polynomial $f$  in $\,x,y,z,w\,$ can be interpreted as a symmetric 
tensor of format $4 \times 4 \times 4$. A typical tensor has complex rank $5$,
but its real rank becomes $6$ as one crosses the
{\em real rank boundary}. This is a
hypersurface of degree $40$  in $\PP^{19}$ studied by
Micha\l ek and Moon \cite[Proposition 3.4]{MM}.

\begin{question}
What can be said about the tropicalization of the Micha\l ek--Moon 
hypersurface?
\end{question}

Seigal \cite[Proposition 2.6]{Sei} identifies the Hessian discriminant as the locus
where the complex rank of cubics $f$ jumps from $5$ to $6$. She points out that 
this discriminant has degree $ \leq 120$ and is invariant
under the action of ${\rm PGL}(4,K)$.

\begin{question}
How to write the Hessian discriminant in terms of Salmon's six
fundamental invariants of the cubic surface?
\end{question}

The  {\em eigenpoints} of $f$ are the fixed points of the
gradient map $\nabla f : \PP^3 \dashrightarrow \PP^3$.
This was studied by Abo et al.~in~\cite{ASS}.
For generic cubics $f$, there are $15$ eigenpoints.
They form the {\em eigenconfiguration} of the $4 {\times} 4 {\times} 4$ tensor~$f$.

\begin{question}
Which configurations of $ 15$
points in $\PP^3$ arise as eigenpoints of a cubic surface?
\end{question}

The {\em eigendiscriminant}  
is a hypersurface of degree $96$ in $\PP^{19}$. This object and its degree are
studied in \cite[Section 4]{ASS}.
It represents cubic surfaces that possess an eigenpoint of multiplicity $\geq 2$.
This hypersuface deserves further study.

\begin{question}
Does there exists a compact determinantal formula for the eigendiscriminant
of the cubic surface?
\end{question}

We learned from
Bruin and Sert\"oz \cite{BrSe}  that there are $255$ Cayley symmetroids
containing a general complete intersection $(2,3)$-curve in $\PP^3$, one for each $2$-torsion
point on the Jacobian of such a genus $4$ curve. This is less than the number $305$ of 
Cayley symmetroids found
in a general $\PP^4$ in $\PP^{19}$; cf.~Question~9.

\begin{question}
What explains the drop from $305$ to $255$ when we count
Cayley symmetroids that lie in the special $4$-plane
in $\PP^{19}$ of all cubic surfaces containing a
given space sextic?
\end{question}

The following question paraphrases Problem 5.4 in \cite{SX}.
 It was studied by Bernal et al.~\cite{BCDFM},
but the authors of that article left it largely unresolved.

\begin{question}
What are all toric degenerations of Cox rings of cubic surfaces?
\end{question}

The following question paraphrases Conjecture 5.3 in \cite{RSS}.

\begin{question}
Can we identify a tropical basis for the universal Cox ideal of cubic surfaces?
\end{question}

In tropical geometry,
it is a big challenge  to relate 
intrinsic Del Pezzo geometry to  embedded geometry in $\PP^3$.
This is reminiscent from the curve case.

\begin{question}
There are two generic types of tropical del Pezzo surfaces of degree $3$,
characterized by the  tree arrangements in  \cite[Figures 4 and 5]{RSS}.
Can we identify cubics $f$ that realize these two types by looking at the valuations of
the six invariants in Question 2?
\end{question}

The lines in $\PP^3$ are points
$p = (p_{12}:p_{13}:p_{14}:p_{23}:p_{24}:p_{34})$
in the Grassmannian ${\rm Gr}(2,4) \subset \PP^5$.
The {\em universal Fano variety} in
$\PP^{19} \times \PP^5$ parametrizes lines on cubic surfaces.
Its ideal is generated modulo the Pl\"ucker quadric by $20$
polynomials of degree $(3,1)$ in $(p,c)$.
These are derived in  \cite[Section 6]{JPS}.

\begin{question}
Can we find an explicit tropical basis for
universal Fano variety?
\end{question}

A real cubic surface in $\PP^3_\RR$ has either one or two connected components.
In the latter case, the cubic is {\em hyperbolic}. It bounds a convex body
that is of interest in optimization. For some background on real cubics we refer to \cite{PBT}.

\begin{question}
Can we find a semialgebraic description for the set of
smooth hyperbolic cubics in $\PP^{19}_\RR$?
How to express this case distinction in terms of the six fundamental invariants?
\end{question}

Every cubic $f$ whose Hessian discriminant (in Question 15) is non-zero
has a unique representation as a sum of five third powers of linear forms,
$\, f \,\,  = \,\, \ell_1^3 \,+ \, \ell_2^3  \,+ \, \ell_3^3  \,+ \, \ell_4^3  \,+ \, \ell_5^3 $.
This is Sylvester's Pentahedral Theorem.
Salmon \cite{Sal} uses this to write the invariants.
 
\begin{question}
Can we find explicit linear forms $\ell_i \in \ZZ[x,y,z,w]$
such that the $2$-adic tropicalization of $V(f)$ is  tropically smooth.
Which unimodular triangulations of the triple tetrahedron $3 \Delta_3$ arise?
\end{question}

If we project a cubic surface $V(f)$ from a general point $p$ on that surface,
then we get a double-covering of $\PP^2$ branched along a quartic curve.
The $28$  bitangents of that curve are the images of the $27$ lines on $V(f)$ plus
one more line which is the exceptional divisor over~$p$.

\begin{question}
Can this correspondence from $27$ to $28$ be understood in tropical geometry? In particular, 
can we see the seven $4$-tuples of tropical bitangents already in  the tropical cubic surface
${\rm Trop}(V(f))$?
\end{question}

The seven $4$-tuples of bitangents of a tropical
plane quartic are studied by Chan and Jiradilok in  \cite{CJ}.
For cubic surfaces,
Brundu and Logar \cite{BL} offer a computational study of $f$
 via the following alternative normal form:
$$ \begin{matrix}
f & = &
a_1 (2 x^2y - 2 xy^2 + xz^2 - xzw - yw^2 + yzw) \,+\,
a_2 (x-w)(xz+yw) \,+ \\ & &
a_3 (z+w)(yw-xz)\,+\,
a_4 (y-z)(xz+yw) \,+\,
a_5(x-y)(yw-xz). \,\,\end{matrix} $$
This amounts to fixing an {\em L-set}, i.e.~a special configuration of five lines
on $V(f)$.

\begin{question}
How to compute the Brundu--Logar normal form in practice?
Can we write  $a_1,a_2,a_3,a_4,a_5$ as rational functions in
$c_1,c_2,\ldots,c_{20}$? What does this tell us tropically?
\end{question}

Here is another interesting normal form, called the
{\em Cayley--Salmon form} of $f$  by Dolgachev \cite{Dol}.
A general cubic surface has $120$ distinct representations 
\begin{equation}
\label{eq:buckley} f \,\, \, = \,\,\, \ell_1 \ell_2 \ell_3 \, + \, m_1 m_2 m_3, 
\end{equation}
where the $\ell_i$ and $m_j$ are linear forms. 
We learned this  from Buckley and Ko\'sir \cite{BK}
who derived it
from the classical construction of {\em Steiner sets}.
 Steiner called the two triples in (\ref{eq:buckley}) a
 {\em Triederpaar}, and he showed that there are $120$ such pairs.

\begin{question}
How to compute the $120$ Cayley--Salmon representations (\ref{eq:buckley})
of a given cubic surface in practice?
Can this be tropicalized?
\end{question}

\section{History}

This section offers a historical introduction to the cubic surface
and its $27$ lines. It is aimed at students and other non-experts.
One goal is to introduce accessible sources for further study of
this beautiful theme in classical algebraic geometry.
Those who seek a textbook introduction are referred to
 the books by Dolgachev \cite[Chapter 9]{Dol}
and Reid \cite[Chapter 7]{Rei}. 
In applications, one considers cubic surfaces over the real numbers,
and for those we refer to
the book by Segre \cite{BSeg}

Arthur Cayley first showed that there are $27$ lines on a general complex cubic surface, and then Salmon showed that there are $27$ lines on every smooth complex cubic surface,  both in the middle of the 19th century.  Their results are remarkable given that nonsingular cubic surfaces are not all projectively equivalent.  In fact, by a simple parameter count, there is a 
$4$-dimensional family of isomorphism classes of cubic surfaces.   Salmon found six fundamental invariants for cubic surfaces.  Their degrees are $8$, $16$, $24$, $32$, $40$ and $100$, and the square of the last one is a polynomial in the other five.  These invariants are homogeneous polynomials in the $20$ coefficients $c_i$ of the cubic form $f$.
They parameterize naturally the four-dimensional family of projective equivalence classes of cubic surfaces.  
For a concrete example, we computed the full
expansion of Salmon's invariant of degree $8$.
It is a sum of $7261$ terms, which looks like this:
$$ \begin{matrix}
192 c_1^2 c_{10}^3 c_{13}^3
-864 c_1^2 c_{10}^3 c_{13}^{\phantom{1}} c_{16}^{\phantom{1}} c_3^{\phantom{1}}
+2592 c_1^2 c_{10}^3 c_ 3^2 c_4^{\phantom{1}}
+1728 c_ 1^2 c_{10}^2 c_{12}^2 c_{13}^{\phantom{1}} c_4^{\phantom{1}}  \smallskip \ \\
-648 c_ 1^2 c_{10}^2 c_{12}^2 c_{16}^2
-288 c_ 1^2 c_{10}^2 c_{12}^{\phantom{1}} c_{13}^2 c_{15}^{\phantom{1}}
+432 c_1^2 c_{10}^2 c_{12}^{\phantom{1}} c_{13}^{\phantom{1}} c_{16}^{\phantom{1}} c_{20}^{\phantom{1}}
+ \,\, \cdots \smallskip \\
 \,\cdots \,\,
+1728 c_3^{\phantom{1}} c_4^2 c_5^2 c_6^{\phantom{1}} c_9^2
+1728 c_3^{\phantom{1}} c_4^2 c_6^2 c_8^2 c_9^{\phantom{1}}
+384 c_3^{\phantom{1}} c_4^{\phantom{1}} c_7^3 c_9^3
+192 c_4^2 c_6^3 c_9^3.
\end{matrix}
$$

A crucial ingredient in the investigations of Cayley and Salmon was the use of  normal forms for $f$.
  The two most important ones were the Cayley-Salmon form and the Sylvester form.  The Cayley-Salmon form, i.e.~general $f$ is the sum of two triple products of linear forms, gives a nice access to the combinatorial structure of the 27 lines.  Steiner \cite{Ste} called the pair of triples of planes defined by the triple products a {\em Triederpaar}, and showed that a general cubic surface has 120 such Triederpaare.
Sylvester's pentahedral form, i.e.~a 
general $f$ is the sum of five cubes of linear forms, 
is the key to Salmon's computation of invariants.

Schl\" afli exploited the combinatorics of the $27$ lines in \cite {Sch1}.  In particular, he described the $36$ double-sixes formed by the $27$ lines. 
 As was shown later by C. Segre in \cite {CSeg} (see also \cite{Dol2}), the double-sixes correspond one to one  to the equivalence classes of linear determinantal presentation of $f$. 
  Subsequently Schl\" afli, in \cite {Sch2}, classified
 irreducible cubic surfaces according to singularities into $22$ different types, and the number of real lines on real smooth cubic surfaces. A double-six of lines determines a pair of birational morphisms of the cubic surface to the 
plane $\PP^2$. Each morphism contracts one set of six lines in the double-six to points
in $\PP^2$, and it maps the other set to conics in $\PP^2$.  This correspondence, 
showing that a smooth cubic surface is isomorphic to 
the projective plane blown up in six points, goes back to Clebsch \cite {Cle}.

The incidence graph of the $27$ lines has a large  automorphism group, identified with the 
Weyl group $E_6$ by Jordan \cite {Jor}.  The automorphism group of the surface is naturally a subgroup of $E_6$.  
The general nonsingular cubic surface has however no nontrivial automorphisms.  
The complex surfaces with nontrivial automorphisms are the ones with at least one Eckardt point, see  \cite {Eck}.
 These are points contained in three of the lines.  A surface with an 
Eckardt point admits a projective involution that fixes the three lines pointwise.  

In his paper \cite{Gei} in the first volume of Mathematische Annalen,
Geiser  noted the relation between the $27$ lines on a cubic surface and 
$28$ bitangents to a plane quartic curve.  The curve is the branch locus of the projection of the surface away
  from a general point on the surface.  
The birational automorphism defined by interchanging the sheets of this double cover is called the {\em Geiser involution}.  

The early results by Sch\" afli on real cubic surfaces inspired the
 production of models of cubic surfaces with up to $27$ real lines.
The culmination of the work of this period is found in B.~Segre's book \cite {BSeg}.   
This book is a comprehensive investigation into the five different kinds of real smooth cubic surfaces. 
To this date, the book \cite{BSeg} is the best source
 on cubic surfaces over the real numbers. 
 
Well into the $20$th century, algebraic geometry 
textbooks featured the smooth cubic surfaces. 
Baker's book {\em Principles of Geometry} \cite {Bak} is an example.
 This changed in the second half of the $20$th century.
Sheaves and schemes now provided better tools for moduli spaces, like Hilbert schemes and 
geometric invariant theory (GIT).  Dolgachev, van Geemen and Kondo used GIT to describe a $4$-dimensional proper moduli space of $k$-nodal cubic surfaces~\cite{DGK} with $k=0$ or $k=1,...,4$.  Allcock, Carlson and Toledo \cite{ACT1}
constructed a GIT model of the moduli space as quotient of a complex $4$-dimensional ball by a reflection group.  

The textbooks on abstract algebraic geometry by
Hartshorne \cite{Har} and Griffiths and Harris \cite {GH} 
introduce smooth cubic surfaces as embeddings of the plane blown up in six points.  
 Miles Reid's undergraduate textbook  \cite {Rei} gives a more elementary introduction to smooth cubic surfaces and its $27$ lines, without reference to the isomorphism with the plane blown up in six points.
 A comprehensive treatment starting with the $27$ lines and the group $E_6$, covering 
 determinantal representations, power sum presentations and automorphisms, is given 
 by Dolgachev \cite{Dol}.  It includes also a very nice historical overview on cubic surfaces.  
  A nice overview of real cubic surfaces from a modern minimal program perspective is found 
  in Koll\'ar's lecture notes \cite{Kol} on real algebraic surfaces.

Tropical geometry started at the beginning of the $21$st century and it offers
a new perspective on classical algebraic geometry.
  The book by Maclagan and Sturmfels  offers an introduction
which includes tropical surfaces in $3$-space in \cite[Section 4.5]{MS}.  These are the natural tropicalizations of affine surfaces in the complement of the coordinate hyperplanes in $\PP^3$. 
The relationship between complex and tropical geometry is an active
area of research.
  Smooth tropical cubic surfaces may have infinitely many lines \cite{PV, Vig}.
  Recent work in \cite{CD, JPS, RSS}
      concerns different representations of the cubic surface, aimed
      at revealing all the  tropical lines in the tropicalization.  
   The tropicalization of the complement of a triangle in a cubic surface has been considered 
   by Gross et al.~\cite{GHKS} to explain mirror symmetry phenomena of Calabi-Yau varieties.

\bigskip \bigskip

\section{Progress}

Considerable progress on the $27$ questions was made
 between January 2019 and September 2019. 
 This resulted in $14$ articles, written for the special volume of
{\em Le Matematiche}. Here we introduce these papers, with
emphasis on how they address the
prompts in Section 2.  Naturally, many questions remain
unanswered. We conclude with a list of ten open problems, 
extracted from the $27$ questions. We 
consider these to be especially interesting for further research.

The first three papers concern
the group action of ${\rm PGL}(4)$ on the projective space
$\PP^{19}$ of cubic surfaces.
Brustenga, Timme and Weinstein \cite{BTW} use
numerical algebraic geometry to answer Question 1.
The general orbit in $\PP^{19}$ is a $15$-dimensional variety 
of degree $96120$. Cazzador and Skauli \cite{CS} develop the
intersection-theoretic approach to this problem. It rests on
the solution by Aluffi and Faber for
the same problem concerning ${\rm PGL}(3)$-orbits of plane curves.
Elsenhans and Jahnel \cite{EJ} resolve the first part of Question 2, by
giving a practical algorithm for evaluating invariants and covariants of cubic surfaces.
Based on the classical method of transvections, it is implemented in {\tt Magma}.

We next come to the discriminant of the cubic surface.
Kastner and L\"owe  \cite{KL} present a computational solution
to Question 3: the Newton polytope has 
$166104$ vertices. The paper by Bunnett and Keneshlou \cite{BuKe} addresses
Questions 4, 6 and 7. The answer to Question 4 is ``no'' since
there is no rank $1$ Ulrich sheaf on the Veronese surface in $\PP^9$. 
The authors determine a rank $2$ Ulrich sheaf, they construct
a Pfaffian representation of size $16 \times 16$ for the discriminant,
and they examine the rank strata of Nansen's $20 \times 20$ matrix.
Keneshlou  \cite{Ken} identifies the
    singular locus of the Eckardt hypersurface. This solves Question 13.

Normal forms of cubic surfaces are important for many applications.
Panizzut, Sert\"oz and Sturmfels \cite{PSS} introduce 
a new normal form  which works well for tropical geometry.
They also answer most of Question 11.
Donten-Bury, G\"orlach and Wrobel \cite{DGW}
resolve Question 19 by describing a
classification of all toric degenerations of 
cubic surfaces. Their approach rests on Khovanskii bases of
Cox rings.  Hahn, Lamboglia and Vargas \cite{HLV} address
Question 27. They describe two methods for computing
 the $120$ Cayley-Salmon equations.
 
 Cubic surfaces correspond to symmetric tensors of format $4 \times 4 \times 4$.
Seigal and Sukarto \cite{SS} investigate how their tensor rank
is reflected in the singularity structure of the surface.
An important player in their paper is the Hessian discriminant. 
This invariant is the object studied by
Dinu and Seynnaeve \cite{DS}, who present the answer to
 Question 15. The article by \c{C}elik, Galuppi, Kulkarni and Sorea \cite{CGKS}
studies the spectral theory of symmetric $4 \times 4 \times 4$ tensors. Their
parametrization of the eigenpoints furnishes a
partial answer to Question~16.

Two articles examine cubic surfaces via
tropical geometry. Brandt and Geiger \cite{BrGe} 
give a partial answer to Question 10. They develop
a tropical theory of binodal cubics
through $17$ given points in $\PP^3$.
In that setting, the classical count of $280$ drops to $214$.
Geiger \cite{Geig} studies the combinatorics of  lines on tropical cubic surfaces.
She proves that the Brundu-Logar normal form 
gives surfaces that are not tropically smooth, thus 
answering part two of Question~27.

We also note that an observant referee proposed an answer for Question 18.
The following explanation for the drop was given. Consider the
special $4$-space of cubic surfaces that contain a  $(2,3)$-curve $C$.
This $4$-space contains a $\PP^3$ of reducible symmetroids, all with 
the unique quadric surface $Q$ that contains $C$, as a component. 
This $\PP^3$ is an excess intersection that counts as $50$ in the $305$.

\smallskip

The $14$ articles offer considerable new insights 
on cubic surfaces. However,
some of the harder questions still remain unanswered.
We conclude with a  reprise of Section 2, in the form of
a list of  {\bf ten open problems} on cubic surfaces.

\begin{enumerate}
\item Determine the prime ideal of the generic ${\rm PGL}(4)$ orbit on $\PP^{19}$. \hfill (Q1)
\vspace{-0.1in}
\item How to evaluate the six fundamental invariants for a tropical cubic?~(Q2) \vspace{-0.1in}
\item Identify all rank varieties for the 
discriminant matrices in \cite{BuKe}. \hfill (Q5-6)
\vspace{-0.1in}
\item Determine the prime ideals of the varieties of $k$-nodal cubics in $\PP^{19}\!$.~(Q7)
\vspace{-0.1in}
\item Find $19-k$ points whose interpolating $k$-nodal surfaces
 are all real.~(Q8-9)
 \vspace{-0.26in}
\item Prove \cite[Conjecture 4.1]{PSS}: 
{\em Smooth tropical cubics have $27$ lines}. \hfill (Q11)
\vspace{-0.1in}
\item Which cubics are tropically smooth after a coordinate change? \hfill (Q12)
\vspace{-0.1in}
\item Find determinantal or pfaffian formulas for eigendiscriminants. \hfill (Q17)
\vspace{-0.1in}
\item What is the correct tropicalization of Sylvester's Pentahedral Form? \hfill (Q24)
\vspace{-0.1in}
\item Relate the $27$ lines to the $28$ bitangents in the tropical 
setting. \hfill (Q25)
\end{enumerate}

These are our favorites among problems extracted from the $27$ questions.
They underscore our belief that, even two centuries after Cayley and Salmon,
the investigation of cubic surfaces will continue to be an active area of research.

\begin{small}

\end{small}


\begin{thebibliography}{10}

\setlength{\itemsep}{-0.2mm}

\bibitem{ASS} H.~Abo, A.~Seigal and B.~Sturmfels:
{\em Eigenconfigurations of tensors}, Algebraic and Geometric Methods in Discrete Mathematics, Contemporary Mathematics 685, Amer. Math. Soc., Providence, RI (2017) 1--25.

\bibitem{ACT1} D.~Allcock, J. A.~Carlson, D.~Toledo: {\em The complex hyperbolic geometry of the moduli space of cubic surfaces}, J. Algebraic Geom. {\bf 11} (2002), 659-724.

\bibitem{ACT2} D.~Allcock, J.A.~Carlson, D.~Toledo: {\em Hyperbolic geometry and
 moduli of real cubic surfaces}. Ann. Sci. \'Ec. Norm. Sup\'er. (4) 43 (2010), no. 1, 69--115.

\bibitem{Bak} H.~Baker:  {\em Principles of Geometry}, vols.1-6., Cambridge University Press, 1922.

\bibitem{Bek} N.D.~Beklemishev: {\em Invariants of cubic forms of four variables},
Vestnik Moscow Univ.~Ser.~I Mat.~Mekh. {\bf 2} (1982) 42--49.

\bibitem{BCDFM} M.~Bernal, D.~Corey, M.~Donton-Bury, N.~Fujita and G.~Merz:
{\em Khovanskii bases of Cox-Nagata rings and tropical geometry},
Combinatorial Algebraic Geometry, 133--157,
Fields Inst.~Communications {\bf 80}, Springer, New York.

 \bibitem{BT}
 P.~Breiding and S.~Timme: {\em
HomotopyContinuation.jl: A package for homotopy continuation in Julia},
 Lecture Notes in Computer Science {\bf 10931},  458-465, 2018.

\bibitem{BrSe} N.~Bruin and E.~Sert\"oz:
{\em Prym varieties of genus four curves},
    Transactions of the American Mathematical Society (2020), 
{\tt arXiv:1808.07881}.

\bibitem{BL} M.~Brundu and A.~Logar:
{\em Parametrization of the orbits of cubic surfaces},
Transformation Groups {\bf 3} (1998) 209--239.

\bibitem{BK} A.~Buckley and T.~Ko\'sir:
{\em Determinantal representations of smooth cubic surfaces} {\bf 125} (2007) 115--140.

\bibitem{CJ} M.~Chan and P.~Jiradilok: {\em Theta characteristics of tropical $K_4$-curves},
Combinatorial Algebraic Geometry, 65--86,
Fields Inst.~Comm.~{\bf 80}, Springer, New York.

\bibitem{Cle} A.~Clebsch: {\em Die Geometrie auf den Fl\"achen dritter Ordnung}, J. Reine Angew. Math.~{\bf 65} (1866) 359-380.

\bibitem{Cre} L.~Cremona: {\em M\'emoire de g\'eom\'etrie pure sur les surfaces du troisi\'eme ordre}, Journal des Math.~pures et appl.~{\bf 68} (1868)

\bibitem{CD} M.A.~Cueto, A.~Deopurkar: {\em Anticanonical tropical cubic del Pezzos contain exactly 27 lines}, {\tt arXiv:1906.08196}

\bibitem{DK} H.~Derksen and G.~Kemper: {\em Computational Invariant Theory},
Encyclopedia of Mathematical Sciences, vol 130, Springer-Verlag, Berlin 2002.    

\bibitem{Dol2} I. V. ~Dolgachev: {\em Luigi Cremona and cubic surfaces}. In Luigi Cremona (1830-1903). Convegno
di Studi matematici, Istituto Lombardo, Accademia di Scienze e Lettere, Milano (2005), 55-70

\bibitem{Dol} I.V.~Dolgachev:  {\em Classical Algebraic Geometry, A Modern View}, Cambridge University Press, 2012.

\bibitem{DGK} I.V.~ Dolgachev, B. ~van Geemen, S. ~Kondo: {\em A complex ball uniformization of the moduli space of cubic surfaces via periods of K3 surfaces}, J. Reine Angew. Math. {\bf 588} (2005), 99--148.

\bibitem{Eck} F.E.~Eckardt: {\em Uber diejenigen Fl\"achen dritten Grades, auf denen sich drei gerade Linien in einem Punkte schneiden}, Mathematische Annalen  {\bf 10} (1876) 227-272

\bibitem{edge} W.L.~Edge: {\em The discriminant of a cubic surface},
Proceedings of the Royal Irish Academy Ser.~A {\bf 80A} (1980) 75--78.

\bibitem{Gei} C.~Geiser: {\em \"Uber die Doppeltangenten einer 
ebenen Curve vierten Grades}, Mathematische Annalen {\bf 1} (1860) 129-138.

\bibitem{GKZ} I.~M.~Gel'fand, M.~Kapranov and A.~Zelevinsky:
{\em Discriminants, Resultants, and Multidimensional Determinants},
Birkh\"auser, Boston, MA, 1994.

\bibitem{GH} P.~Griffiths and J.~Harris: {\em Principles of Algebraic Geometry}, 
Wiley \& Sons, 1978.

\bibitem{GHKS} M.~Gross, P.~Hacking, S.~Keel, B.~Siebert: {\em The mirror of the cubic surface}, {\tt arXiv:1910.08427}

\bibitem{HKT} P.~Hacking, S.~Keel and J.~Tevelev:
{\em Stable pair, tropical, and log canonical compactifications of moduli spaces of del Pezzo 
surfaces}, Invent.~Math.~{\bf 178} (2009) 173--227.

\bibitem{Har} R.~Hartshorne:  {\em Algebraic Geoometry}, GTM {\bf 52}, Springer (1977)

\bibitem{Jor} C.~Jordan: {\em Trait\'e des substitutions et \'equations alg\'ebriques}, Paris, Gauthier- Villars, 1870.

\bibitem{JPS} M.~Joswig, M.~Panizzut and B.~Sturmfels:
{\em The Schl\"afli fan}, {\tt arXiv:1905.11951}.

\bibitem{Kaz} B.Ya.~Kazarnovskii: {\em Newton polyhedra
and the B\'ezout formula for matrix-valued functions of finite-dimensional
representations}, Functional Analysis and its Applications {\bf 21} (1987) 319--321.



\bibitem{Kol} J.~Koll\'ar:
{\em Real algebraic surfaces}, Lecture notes of Trento summer school, September 1997,
{\tt arXiv:alg-geom/9712003}.


\bibitem{MS} D.~Maclagan and B.~Sturmfels:
{\em Introduction to Tropical Geometry},
Graduate Studies in Mathematics, Vol 161, American Mathematical Society, 2015. 

\bibitem{MM} M.~Micha\l ek and H.~Moon:
{\em Spaces of sums of powers and real rank boundaries},
Beitr\"age zur Algebra und Geometrie {\bf 59} (2018) 645--663.

\bibitem{Nan} E.~J.~Nansen: {\em On the eliminant of a set of quadrics, ternary or 
quaternary}, Proceedings of the Royal Society of Edinburgh {\bf 22} (1899) 353--358.

\bibitem{PV} M.~Panizzut and M.D.~Vigeland: {\em Tropical Lines on Cubic Surfaces}, 
{\tt arXiv:0708.3847} (revised 2019).

\bibitem{PBT} I.~Polo-Blanco and J.~Top:
{\em Explicit real cubic surfaces},
Canadian Mathematical Bulletin {\bf 11} (2008) 125--133.

\bibitem{Rei} M.~Reid: {\em Undergraduate Algebraic Geometry}, LMS Student texts {\bf 12}, Cambridge University Press, 1988.

\bibitem{RSS} Q.~Ren, K.~Shaw and B.~Sturmfels:
{\em Tropicalization of del Pezzo surfaces},
Advances in Mathematics {\bf 300} (2016) 156--189. 

\bibitem{Sal} G.~Salmon: {\em On quaternary cubics},
Philosophical Transactions of the Royal Society,
{\bf 150} (1861) 229--239.

\bibitem{Sch1} L.~Schl\"afli: {\em An attempt to determine the twenty-seven lines upon a surface of the third order and to divide such surfaces into species in reference to the reality of the lines upon the surface}, Quart. J. Math. {\bf 2} (1858)
 55--65, 110--121.

\bibitem{Sch2} L.~Schl\"afli: {\em On the distribution of surfaces of the third order into species, in reference to the absence or presence of singular points, and the reality of their lines}, Phil. Trans. of Roy. Soc. London {\bf 6} (1863) 201--241.

\bibitem{BSeg} B.~Segre: {\em The Non-singular Cubic Surfaces}, 
Oxford University Press, 1942

\bibitem{CSeg} C.~Segre: {\em Sur la g\'en\'eration projective des surfaces cubiques}, 
Archiv der Math.~und Phys.~{\bf 10} (1906) 209--215 

\bibitem{Sei} A.~Seigal: {\em Ranks and symmetric ranks of cubic surfaces},
{\tt arXiv:1801.05377}.

\bibitem{Ste} J.~Steiner: {\em \"Uber die Fl\"achen dritten Grades},
 Journ. f\"ur reiner und angew. Math., {\bf 53} (1856)e 133--141

\bibitem{SX} B.~Sturmfels and Z.~Xu: {\em Sagbi bases of Cox-Nagata rings},
  Journal of the European Mathematical Society {\bf 12} (2010) 429--459. 

\bibitem{Vai} I.~Vainsencher: {\em Hypersurfaces with up to six double points},
Communications in Algebra {\bf 31} (2003) 4107--4129.

\bibitem{Vig} M.D.~Vigeland: {\em Smooth tropical surfaces with infinitely 
many tropical lines}, Ark. Mat. 48 (2010),  {\bf 1}, 177-206.

\medskip

\bibitem{BrGe} M.~Brandt and A.~Geiger: {\em A tropical count of binodal cubic surfaces}, {\tt arXiv:1909.09105}.

\bibitem{BTW} L.~Brustenga, S.~Timme, M.~Weinstein: {\em 96120: The degree of the linear orbit of a cubic surface}, {\tt arXiv:1909.06620}.

\bibitem{BuKe} D.~Bunnett and H.~Keneshlou: {\em Determinantal representations of the cubic discriminant}.
{\tt arXiv:1909.05579}.

\bibitem{CS} E.~Cazzador and B.~Skauli: {\em Towards the degree of the $PGL(4)$-orbit of a cubic surface}, preprint.

\bibitem{CGKS} T.~\c{C}elik, F.~Galuppi, A.~Kulkarni and M.-S.~Sorea: {\em On the eigenpoints of cubic surfaces},{\tt arXiv:1909.06261}.

\bibitem{DS} R.~Dinu and T.~Seynnaeve: {\em The Hessian discriminant}, {\tt arXiv:1909.06681}.

\bibitem{DGW} M.~Donten-Bury, P.~G\"orlach and M.~Wrobel: 
  {\em Towards classifying toric degenerations of cubic surfaces}, {\tt arXiv:1909.06690}.

\bibitem{EJ} A.-S.~Elsenhans and J.~Jahnel: {\em Computing invariants of cubic surfaces}, {\tt arXiv:1909.00497}.

\bibitem{Geig} A.~Geiger: {\em On realizability of lines on tropical cubic surfaces and the Brundu-Logar normal form}, {\tt arXiv:1909.09391}.

\bibitem{HLV} M.~Hahn, S.~Lamboglia and A.~Vargas: {\em A short note on Cayley-Salmon equations}, preprint.

\bibitem{KL} K.~Kastner and R.~L\"owe: {\em The Newton polytope of the discriminant of a quaternary cubic form}, {\tt arXiv:1909.08910}.

\bibitem{Ken} H.~Keneshlou:  {\em Cubic surfaces on the singular locus of the Eckardt hypersurface}, {\tt arXiv:1909.05554}.

\bibitem{PSS} M.~Panizzut, E.~Sert\"oz and B.~Sturmfels: {\em An octanomial model for cubic surfaces}, {\tt arXiv:1908.06106}.

\bibitem{SS} A.~Seigal and E.~Sukarto: {\em Ranks and singularities of cubic surfaces}, {\tt arXiv:1909.12538}.

\end{thebibliography}
\end{document}